\documentclass[11pt]{amsart}
\usepackage{geometry}                % See geometry.pdf to learn the layout options. There are lots.
\usepackage{graphicx}
\usepackage{amssymb}
\usepackage{epstopdf}
\usepackage{lscape}
\usepackage{setspace}
\title{Properties of reciprocity formulas for the Rogers-Ramanujan continued fractions}
\author{Rajeev Kohli}
\date{February 22, 2019}
\begin{document}	
\maketitle
\centerline{Columbia University
%\footnote{Address: 506 Uris Hall, Columbia University, New York, NY 10027. Email: rk35@columbia.edu}
}

\begin{abstract}
Ramanujan recorded four reciprocity formulas for the Roger-Ramanujan continued fractions.  Two reciprocity formulas each are also associated with the Ramanujan--G\"ollnitz--Gordon continued fractions  and a  level-13 analog of the Roger-Ramanujan continued fractions. We show that all eight reciprocity formulas are related to a pair of quadratic equations. The solution to the first equation generalizes the golden ratio and is used to set the value of a coefficient in the second equation; and the solution to the second equation gives a pair of values for a continued fraction. We relate the coefficients of the quadratic equations to important formulas obtained by Ramanujan, examine a pattern in the relation between a continued fraction and its parameters, and use the reciprocity formulas to obtain close approximations for all values of the continued fractions. We highlight patterns in the expressions for certain explicit values of the Rogers-Ramanujan continued fractions by expressing them in terms of the golden ratio. We extend the analysis to reciprocity formulas for Ramanujan's cubic continued fraction and the Ramanujan-Selberg continued fraction.

\keywords{Roger-Ramanujan continued fractions \and  Ramanujan--G\"ollnitz--Gordon continued fractions \and Ramanujan's cubic continued fractions \and Ramanujan-Selberg continued fraction \and Reciprocal formulas}
\end{abstract}

\section{Introduction}
Ramanujan (\cite{rln}, \cite{Ram}) recorded four reciprocity formulas for the Rogers-Ramanujan continued fractions. Two other reciprocity formulas are associated with the Ramanujan--G\"ollnitz--Gordon continued fractions (Bagis \cite{Bag}), and another two with a level-13 analog of the continued fractions (Cooper and Ye \cite{CY}). We show that these eight formulas are obtained by choosing different values of a single parameter, $k$, in the quadratic equation $c^2-kc-1=0$, and then substituting the value of its root $c=(k+\sqrt{k^2+4})/2$ into the equation $(1-x_1x_2)/(x_1+x_2)=c$. The eight parameter values are  $k=\pm 1, \pm 2, \pm 3$ and $\pm 11$. Table \ref{tab1} shows the terms associated with $x_1$ and $x_2$ in a reciprocity formula for each value of $k$. In the table, $R(q)$ and $S(q)$ denote the Rogers-Ramanujan continued fractions; $V(q)$ and $-V(-q)$ denote the Ramanujan--G\"ollnitz--Gordon continued fractions; and $R'(q)$ and $S'(q)$ denote the level-13 analogs of the Rogers-Ramanujan continued fractions. The values of $c$ in the table are $\phi=(\sqrt{5}+1)/2$, $1/\phi=(\sqrt{5}-1)/2$, $\phi^5=(5\sqrt{5}+11)/2$, $1/\phi^5=(5\sqrt{5}-11)/2$, $\nu=\sqrt{2}+1$, $1/\nu=\sqrt{2}-1$, $\rho=(\sqrt{13}+3)/2$ and $1/\rho=(\sqrt{13}-3)/2$. 

\begin{table}[htp]
\caption{Values of $k$ and $c$ associated with the Rogers-Ramanujan continued fractions, the Ramanujan--G\"ollnitz--Gordon continued fractions and their level-13 analog}
\begin{center}
\begin{tabular}{rrrrr}
\hline
$k$&$c$&$x_1$&$x_2$&$\alpha\beta$\\
\hline
$1$&$\phi$&$R\left(e^{-2\alpha}\right)$&$R(e^{-2\beta})$&$\pi^2$\\
$-1$&$1/\phi$&$S(e^{-\alpha})$&$S(e^{-\beta})$&$\pi^2$\\
$11$&$\phi^5$&$R^5\left(e^{-2\alpha}\right)$&$R^5\left(e^{-2\beta}\right)$&$\pi^2/5$\\
$-11$&$1/\phi^5$&$S^5(e^{-\alpha})$&$S^5(e^{-\beta})$&$\pi^2/5$\\
$2$&$\nu$&$V(e^{-\alpha})$&$V(e^{-\beta})$&$\pi^2$\\
$-2$&$1/\nu$&$-V(-e^{-\alpha})$&$-V(-e^{-\beta})$&$\pi^2$\\
$3$&$\rho$&$R'(e^{-2\alpha})$&$R'(e^{-2\beta})$&$\pi^2/13$\\
$-3$&$1/\rho$&$S'(e^{-\alpha})$&$S'(e^{-\beta})$&$\pi^2/13$\\
\hline
\end{tabular}
\end{center}
\label{tab1}
\end{table}
We use the two equations to examine the common properties of the continued fractions. If $x_1=x_2=x$, then $(1-x_1x_2)/(x_1+x_2)=c$ reduces to $(1/x)-x=c$. This form has been used to express important formulas for $x=R(q), R^5(q), V(q)$ and $R'(q)$.  If we set $x_1+x_2=2r$, then $(1-x_1x_2)/(x_1+x_2)=c$  can be written as the quadratic equation $x^2-2rx-(2rc-1)=0$, its two roots corresponding to the pairs of continued fractions shown in Table \ref{tab1}. The value of $r$ lies between $1/2c$ and $\sqrt{c^2+1}-c$. We examine the patterns in the relations between $r$, $x_1$, $x_2$, $k$ and $q$ when $x_1$ and $x_2$ represent the pairs of continued fractions in Table \ref{tab1}.  In each case, as $\alpha$ increases, $r$ first decreases and then increases, attaining a minimum value when $\alpha=\beta$. For example, consider $k=1$, $x_1=R\left(e^{-2\alpha}\right)$ and $x_2=R(e^{-2\beta})$. If $\alpha=0$, $R\left(e^{-2\alpha}\right)=1/\phi$ and  $R(e^{-2\beta})=0$. As $\alpha$ increases from zero to $\alpha=\beta=\pi$, $R\left(e^{-2\alpha}\right)$ decreases and $R(e^{-2\beta})$ increases to $\sqrt{\phi^2+1}-\phi$. The functions $R\left(e^{-2\alpha}\right)$ and $R(e^{-2\beta})$ change very slowly with $\alpha$ when $0<\alpha<\pi/5$. Thereafter, they change almost linearly until $\alpha=\beta=\pi$. Similar observations can be made about the continued fractions in the other pairs of reciprocity formulas. We highlight patterns in the expressions for certain explicit values of the Rogers-Ramanujan continued fractions by expressing them in terms of the golden ratio. Finally, we extend the present analysis to reciprocity formulas for Ramanujan's cubic continued fractions and the Ramanujan-Selberg continued fraction. 

\section{Reciprocity formulas}
The Rogers-Ramanujan continued fractions are
$$R(q)={q^{1/5}\over 1+}{q\over 1+}{q^2\over 1+}{q^3\over 1+}\dots,$$
$$\qquad\qquad\qquad\quad\ \  =q^{1/5}\prod_{j=1}^\infty {(1-q^{5j-1})(1-q^{5j-4})
\over (1-q^{5j-2})(1-q^{5j-3})}$$
$$\qquad\  \ =q^{1/5}\prod_{j=1}^\infty (1-q^j)^{\left({j\over 5}\right)}$$
where $|q|<1$ and $\left({j\over p}\right)$ is the Legendre symbol; and  
$$S(q)=-R(-q)={q^{-1/5}\over 1-}{q\over 1+}{q^2\over 1-}{q^3\over 1+}\dots.$$

\bigskip\noindent
Consider the equation
$${1-x_1x_2\over x_1+x_2}=c.\eqno(1)$$
We write it as 
$x_1x_2+c(x_1+x_2)=1$, then add $c^2$ to both sides to obtain
$$(c+x_1)(c+x_2)=1+c^2.\eqno(2)$$
Let 
$$c^2-kc-1=0.\eqno(3)$$
We substitute its root
$$c={k+\sqrt{k^2+4}\over 2}\eqno(4)$$
into equation (2) to obtain
$$\left({\sqrt{k^2+4}+k\over 2}+x_1\right)\left({\sqrt{k^2+4}+k\over 2}+x_2\right)=\sqrt{k^2+4}\left({\sqrt{k^2+4}+k\over 2}\right).\eqno(5)$$

\bigskip\noindent
1. Let $k=\pm 1$. Then equation (5) becomes 
$$\quad\ \  \left(\phi+x_1\right)\left(\phi+x_2\right)=\left(\phi+{1\over \phi}\right)\phi,\eqno(6)$$
$$\left({1\over \phi}+x_1\right)\left({1\over \phi}+x_2\right)=\left(\phi+{1\over \phi}\right){1\over \phi},\eqno(7)$$
where $\phi=(\sqrt{5}+1)/2$, $1/\phi=(\sqrt{5}-1)/2$ and $\phi+(1/\phi)=\sqrt{5}$. Entry 39 in Chapter XVI of Ramanujan's second notebook \cite{Ram} records the pair of reciprocity formulas in equations (6) and (7), with $x_1=R\left(e^{-2\alpha}\right)$ and $x_2=R(e^{-2\beta})$ in equation (6), and $x_1=S(e^{-\alpha})$ and $x_2=S(e^{-\beta})$ in equation (7), where $\alpha\beta=\pi^2$.

\bigskip\noindent
2. Let $k=\pm 11$. Then equation (5) becomes
$$\quad \ \left(\phi^5+x_1\right)\left(\phi^5+x_2\right)=\left(\phi^5+{1\over \phi^5}\right)\phi^5,\eqno(8)$$
$$\left({1\over \phi^5}+x_1\right)\left({1\over \phi^5}+x_2\right)
=\left(\phi^5+{1\over \phi^5}\right){1\over \phi^5},\eqno(9)$$
where $\phi^5=(5\sqrt{5}+11)/2$, $1/\phi^5=(5\sqrt{5}-11)/2$ and $\phi^5+(1/\phi^5)=5\sqrt{5}$. Entry 3.2.9 in Ramanujan's lost notebook (\cite{rln}, p. 91) records the reciprocity formulas in equations (8) and (9), with $x_1=R^5\left(e^{-2\alpha}\right)$ and $x_2=R^5\left(e^{-2\beta}\right)$ in equation (8), and $x_1=S^5(e^{-\alpha})$ and $x_2=S^5(e^{-\beta})$ in equation (9), where $\alpha\beta=\pi^2/5$.

Since $\phi^5=5\phi+3$ and $1/\phi^5=5\phi-8$, equation (1) implies
$$\left({1-R\left(e^{-2\alpha}\right)R(e^{-2\beta})\over R\left(e^{-2\alpha}\right)+R(e^{-2\beta})}\right)^5=5\left({1-R\left(e^{-2\alpha}\right)R(e^{-2\beta})\over R\left(e^{-2\alpha}\right)+R(e^{-2\beta})}\right)+3={1-R^5(e^{-2\gamma})R^5(e^{-2\delta})\over R^5(e^{-2\gamma})+R^5(e^{-2\delta})}$$
and
$$\left({1-S^5(e^{-\alpha})S^5(e^{-\beta})\over S^5(e^{-\alpha})+S^5(e^{-\beta})}\right)^5=5\left({1-S(e^{-\alpha})S(e^{-\beta})\over S(e^{-\alpha})+S(e^{-\beta})}\right)-8={1-S^5(e^{-\gamma})S^5(e^{-\delta})\over S^5(e^{-\gamma})+S^5(e^{-\delta})}$$
where $\alpha\beta=\pi^2$ and $\gamma\delta=\pi^2/5$.  

\bigskip\noindent
3. Let $k=\pm 2$. Then equation (5) becomes
$$\left(\nu+x_1\right)\left(\nu+x_2\right)=\left(\nu+{1\over \nu}\right)\nu,\eqno(10)$$
$$\left({1\over \nu}+x_1\right)\left({1\over \nu}+x_2\right)=\left(\nu+{1\over \nu}\right){1\over \nu},\eqno(11)$$
where $\nu=\sqrt{2}+1$, $1/\nu=\sqrt{2}-1$ and $\nu+(1/\nu)=2\sqrt{2}$. Ramanujan (\cite{Ram}, Chapter XIX, p. 225) recorded the following continued fraction, which is now known as the Ramanujan--G\"ollnitz--Gordon continued fraction (Chan and Huang \cite{CH}, Baruah and Saikia \cite{BS}): 
$$V(q)={q^{1/2}\over 1+q+}{q^2\over 1+q^3+}{q^4\over 1+q^5+}{q^6\over 1+}\dots$$
$$\qquad\ \  =q^{1/2}\prod_{j=1}^\infty {(1-q^{8j-7})(1-q^{8j-1})\over (1-q^{8j-5})(1-q^{8j-3})}.$$
The reciprocity formulas in equations (10) and (11) are obtained when $x_1=V(e^{-\alpha})$ and $x_2=V(e^{-\beta})$ in equation (10), and $x_1=-V(-e^{-\alpha})$ and $x_2=-V(-e^{-\beta})$ in equation (11), where $\alpha\beta=\pi^2$ (Bagis \cite{Bag},  Corollary 3.1, p. 11). 

\bigskip\noindent
4. Let $k=\pm 3$. Then equation (5) becomes
$$\quad\ \  \left(\rho+x_1\right)\left(\rho+x_2\right)=\left(\rho+{1\over \rho}\right)\rho,\eqno(12)$$
$$\left({1\over \rho}+x_1\right)\left({1\over \rho}+x_2\right)=\left(\rho+{1\over \rho}\right){1\over \rho},\eqno(13)$$
where $\rho=(\sqrt{13}+3)/2$, $1/\rho=(\sqrt{13}-3)/2$ and $\rho+(1/\rho)=\sqrt{13}$.  
In his second notebook, Ramanujan recorded the following level-13 analog of $R(q)$: 
$$R'(q)=q\prod_{j=1}^\infty {(1-q^{13j-12})(1-q^{13j-10})(1-q^{13j-9})(1-q^{13j-4})(1-q^{13j-3})(1-q^{13j-1})\over (1-q^{13j-11})(1-q^{13j-8})(1-q^{13j-7})(1-q^{13j-6})(1-q^{13j-5})(1-q^{13j-2})}$$
$$=q\prod_{j=1}^\infty (1-q^j)^{\left({j\over 13}\right)}.\qquad\qquad\qquad\qquad\qquad\qquad\qquad\qquad\qquad\quad\ \ \ $$
Cooper and Ye \cite{CY} obtained the reciprocity formulas in equations (12) and (13) with $x_1=R'(e^{2\alpha})$ and $x_2=R'(e^{2\beta})$ in equation (12),  and $x_1=S'(e^{2\alpha})$ and $x_2=S'(e^{2\beta})$ in equation (13), where $S'(q)=-R'(-q)$ and $\alpha\beta=\pi^2/13$.

\bigskip\noindent
Andrews and Berndt \cite{rln} observed that the following are among the most fundamental and useful properties of the Rogers-Ramanujan continued fraction: 
$$\quad\ \  {1\over R(q)}-R(q)=1+{f(-q^{1/5})\over q^{1/5}f(-q^{5})}=1+{1\over q^{1/5}}\prod_{j=1}^\infty {(1-q^{j/5})\over (1-q^{5j})}\eqno(14)$$
and
$${1\over R^5(q)}-R^5(q)=11+{f^6(-q)\over qf^6(-q^{5})}=11+{1\over q}\prod_{j=1}^\infty {(1-q^j)^6\over (1-q^{5j})^6},\eqno(15)$$
where 
$$f(-q)=(q; q)_{\infty}=\prod_{j=1}^{\infty} (1-q^j)=\sum_{j=-\infty}^{\infty} (-1)^j q^{j(3j-1)/2}\ \ .$$ 
Ramanujan and subsequently others (e.g., Berndt and Chan \cite{BC}, Soon-Yi \cite{S-Y}) used equations (14) and (15) to evaluate $R(q)$ for specific values of $q$. 

Ramanujan (\cite{Ram}, Chapter XIX, p. 225) recorded the following analogous expressions for $V(q)$: 
$${1\over V(q)}-V(q)={\varphi(q^2)\over q^{1/2}\psi(q^4)}$$
and
$${1\over V(q)}+V(q)={\varphi(q)\over q^{1/2}\psi(q^4)},$$
where 
$$\varphi(q)=\sum_{k=-\infty}^\infty q^{k^2}\qquad $$
and
$$\psi(q)=\sum_{k=0}^\infty q^{k(k+1)/2}\ \ .$$
The corresponding expression for $R'(q)$ is 
$${1\over R'(q)}-R'(q)=3+{f^2(-q)\over qf^2(-13q)}=3+{1\over q}\prod_{j=1}^{\infty} {(1-q^j)^2\over (1-q^{13j})^2}.\eqno(16)$$

The left-hand side in each of these equations is obtained by setting $x_1=x_2=x$ in equation (1) and substituting $x=R(q)$, $x=R^5(q)$, $x=V(q)$ and $x=R'(q)$. In this case 
$${1\over x}-x=2c=k+\sqrt{k^2+4},\eqno(17)$$
which gives
$$x=\sqrt{c^2+1}-c=\sqrt{{k(k+\sqrt{k^2+4})+4\over 2}}-{k+\sqrt{k^2+4}\over 2}.$$

\bigskip\noindent
1. If $x=R(q)$, then 
$$2c-1=k+\sqrt{k^2+4}-1={f(-q^{1/5})\over q^{1/5}f(-q^{5})}\eqno(18)$$
and 
$$k=c-{1\over c},\ \ c={1\over 2}\left(1+{f(-q^{1/5})\over q^{1/5}f(-q^{5})}\right).$$
The value $k=1$ is associated with $q=e^{-2\pi}$, which gives the following relations obtained by Ramanujan:
$${f(-e^{-2\pi/5})\over e^{-2\pi/5}f(-e^{-5\cdot 2\pi})}=2\phi-1=\sqrt{5}$$
and
$$2c=1+{f(-e^{-2\pi/5})\over e^{-2\pi/5}f(-e^{-5\cdot 2\pi})}=1+\sqrt{5}=2\phi.$$

\bigskip\noindent
2. If $x=R^5(q)$, then 
$$k+\sqrt{k^2+4}-11={f^6(-q)\over qf^6(-q^{5})}\eqno(19)$$
and
$$k=c-{1\over c},\ \ c={1\over 2}\left(11+{f^6(-q)\over qf^6(-q^{5})}\right).\eqno(20)$$
The value $k=11$ is associated with $q=e^{-2\pi/\sqrt{5}}$, which gives  
$${f^6(-e^{-2\pi/\sqrt{5}})\over e^{-2\pi/\sqrt{5}} f^6(-e^{-\sqrt{5}\cdot 2\pi})}=5\sqrt{5}$$
and
$$2c={1\over R^5(e^{-2\pi/\sqrt{5}})}-R^5(e^{-2\pi/\sqrt{5}})=11+\sqrt{11^2+4}=11+5\sqrt{5}=2\phi^5.$$

\bigskip\noindent
3. If $x=V(q)$, then 
$$2c={1\over V(q)}-V(q)=k+\sqrt{k^2+4}={\varphi(q^2)\over q^{1/2}\psi(q^4)}\eqno(21)$$
and
$$k=c-{1\over c},\ \ c={\varphi(q^2)\over 2q^{1/2}\psi(q^4)}\ .\eqno(22)$$
The value $k=2$ is associated with $q=e^{-\pi}$, which gives
$$2c={1\over V(e^{-\pi})}-V(e^{-\pi})={\varphi(e^{-2\pi})\over e^{-\pi/2}\psi(e^{-4\pi})}=2\sqrt{2}+2=2\nu.$$

\bigskip\noindent
4. If $x=R'(q)$, then 
$$2c={1\over R'(q)}-R'(q)=k+\sqrt{k^2+4}\eqno(23)$$
and
$$k+\sqrt{k^2+4}-3={1\over q}\prod_{j=1}^{\infty} {(1-q^j)^2\over (1-q^{13j})^2}\ .\eqno(24)$$
The value $k=3$ is associated with $q=e^{-2\pi/\sqrt{13}}$, in which case
$${1\over R'(e^{-2\pi/\sqrt{13}})}-R'(e^{-2\pi/\sqrt{13}})=2\rho.$$

\bigskip\noindent
Equation (3) implies that $c^n=a_n+b_nc$, where $a_1=0$, $b_1=1$, $a_n=b_{n-1}$ and  $b_n=(a_{n-1}+kb_{n-1})$, for all $n\ge 2$. As a result, the left-hand side of equation  (1) is a  linear function of its integer powers.  Thus
$$\left({1-x_1x_2\over x_1+x_2}\right)^n=a_n+b_n\left({1-x_1x_2\over x_1+x_2}\right),\ n\ge 2.\eqno(25)$$
Similarly, $1/c^{n}=A_n+B_nc$, where $A_1=-k$, $B_1=1$, $A_n=B_{n-1}-kA_{n-1}$ and $B_n=A_{n-1}$, for all $n\ge 2$. 
 Thus
$$\left({x_1+x_2\over 1-x_1x_2}\right)^{n}=A_n+B_n\left({1-x_1x_2\over x_1+x_2}\right),\  n\ge 2.\eqno(26)$$
If $x_1=x_2=x$, then
$$\left({1-x^2\over 2x}\right)^n=a_n+b_n \left({1-x^2\over 2x}\right)$$
and
$$\left({2x\over 1-x^2}\right)^{n}=A_n+B_n\left({1-x^2\over 2x}\right)$$
For example, if $n=2$, then 
$$\left({1\over x}-x\right)^2=2^2+2k\left({1\over x}-x\right),$$
which implies that
$${1\over x}-x=2k+{4\over 2k+}{4\over 2k+}{4\over 2k+}{4\over 2k+}\dots.$$
and
$${1\over x^2}+x^2=6+2k\left({1\over x}-x\right)$$
$$\qquad\qquad\qquad \ =6+2k(k+\sqrt{k^2+4}).$$
If $x=R(q)$, then 
$${1\over R^2(q)}+R^2(q)=2+\left(1+{f(-q^{1/5})\over q^{1/5}f(-q^{5})}\right)^2.$$

\section{A second quadratic equation}
Watson \cite{W} showed that if $\alpha\beta=\pi^2$, then 
$R\left(e^{-2\alpha}\right)$ is a root of the quadratic equation 
$${1\over x}-1-x=\left({k_1k_1^{'}K_1^3\over k_2k_2^{'}K_2^3}\right)^{1\over 6}$$
and $R\left(e^{-2\beta}\right)$ is a root of the quadratic equation
$${1\over x}-1-x=\left({k_2k_2^{'}K_2^3\over k_1k_1^{'}K_1^3}\right)^{1\over 6},$$
where $k_1$ and $k_2$ are the moduli of elliptic functions whose quarter-periods satisfy the relations
$$\pi{K_1^{'}\over K_1}={\alpha\over 5},\  \ \pi {K_2^{'}\over K_2}=5\alpha,\ \ 
\pi{K_2\over K_2^{'}}={\beta\over 5},\  \ \pi {K_1\over K_1^{'}}=5\beta.$$

Equation (1) gives a different quadratic equation, the solution to which is a pair of values for the Rogers-Ramanujan continued fractions, the Ramanujan--G\"ollnitz--Gordon continued fractions or the level-13 analog. This equation is obtained as follows. 

Let $x_1+x_2=2r$. Then equation (1) can be written as
$$2rc=1-x_1x_2=1-x_1(2r-x_1)=1-x_2(2r-x_2).$$
That is, $x_1$ and $x_2$ are the roots 
$$x=r\pm \sqrt{r^2+2rc -1}\eqno(27)$$ 
of the quadratic equation 
$$x^2-2rx-(2rc-1)=0.\eqno(28)$$
Observe that any quadratic equation $x^2+Bx+C=0$ can be expressed as equation (28) by setting $r=-B/2$ and $c=(C-1)/B$. Let $2rc\le 1$ so that both the roots in equation (27) are non-negative. If $r^2+2rc -1=0$, then $x_1=x_2=r=\sqrt{c^2+1}-c$.  Otherwise, the larger root has a value between $1/c$ and $\sqrt{c^2+1}-c$, and the smaller a value between $\sqrt{c^2+1}-c$ and zero. 
\begin{figure}[htbp]
\begin{center}
\includegraphics[width=0.75\textwidth]{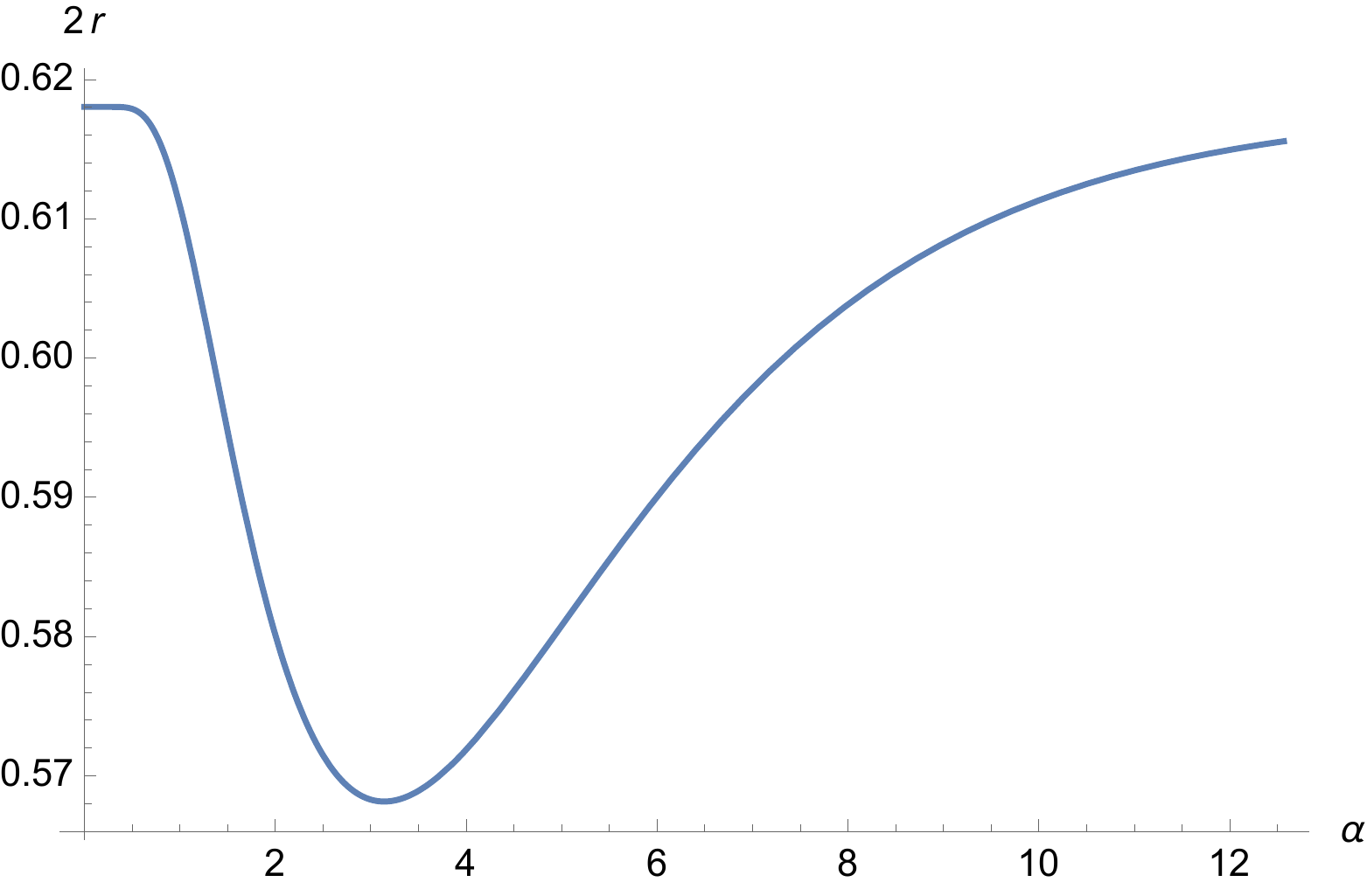}\caption{Value of $2r=R\left(e^{-2\alpha}\right)+R(e^{-2\beta})$ as a function of $\alpha$, where $\alpha\beta=\pi^2$.}
\label{fig1}
\end{center}
\end{figure}
\begin{figure}[htbp]
\begin{center}
\includegraphics[width=0.75\textwidth]{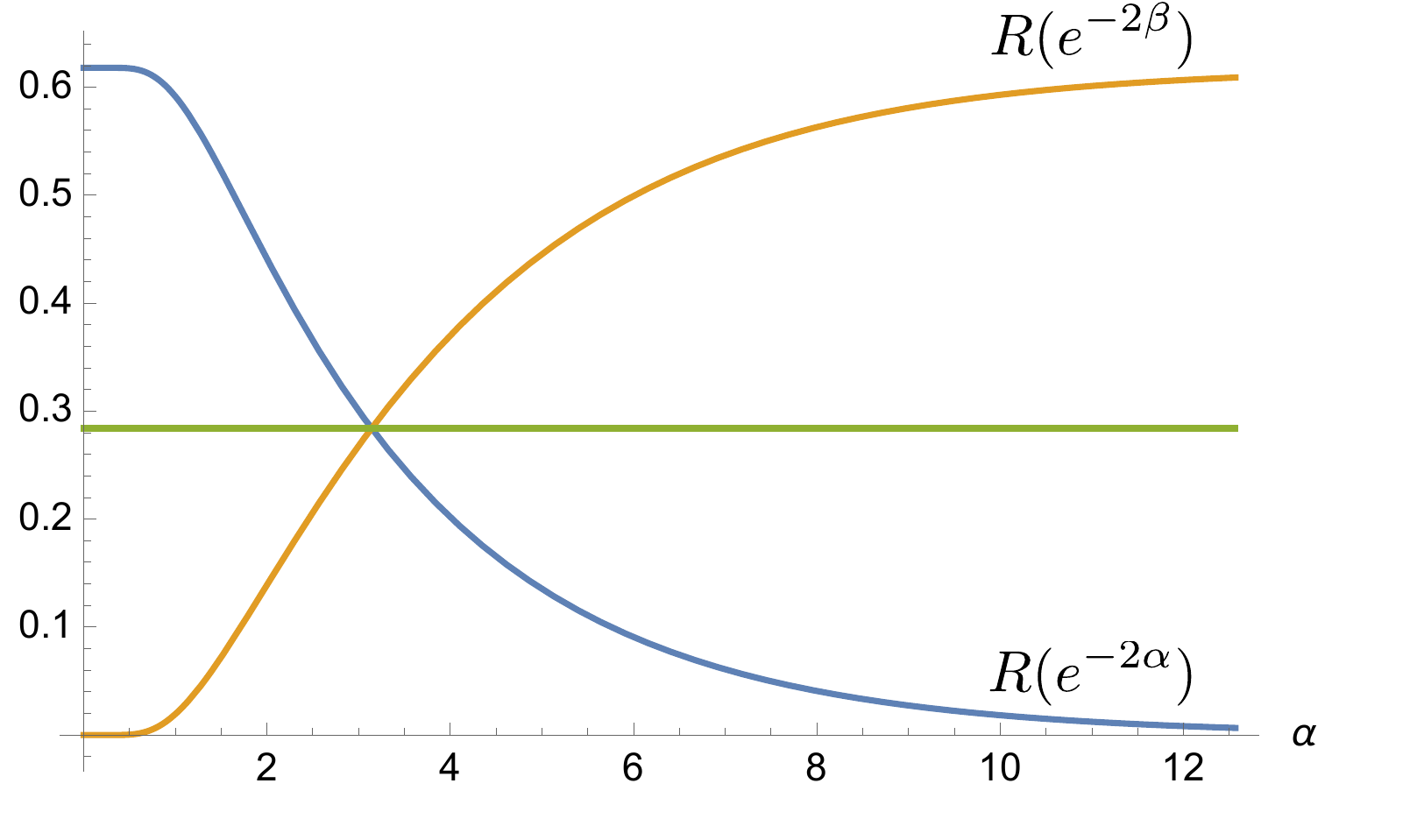}\caption{Values of  $R\left(e^{-2\alpha}\right)$ (blue line) and $R(e^{-2\beta})$ (yellow line) as a function of $\alpha$, where $\alpha\beta=\pi^2$.}
\label{fig2}
\end{center}
\end{figure}
\begin{figure}[htbp]
\begin{center}
\includegraphics[width=0.75\textwidth]{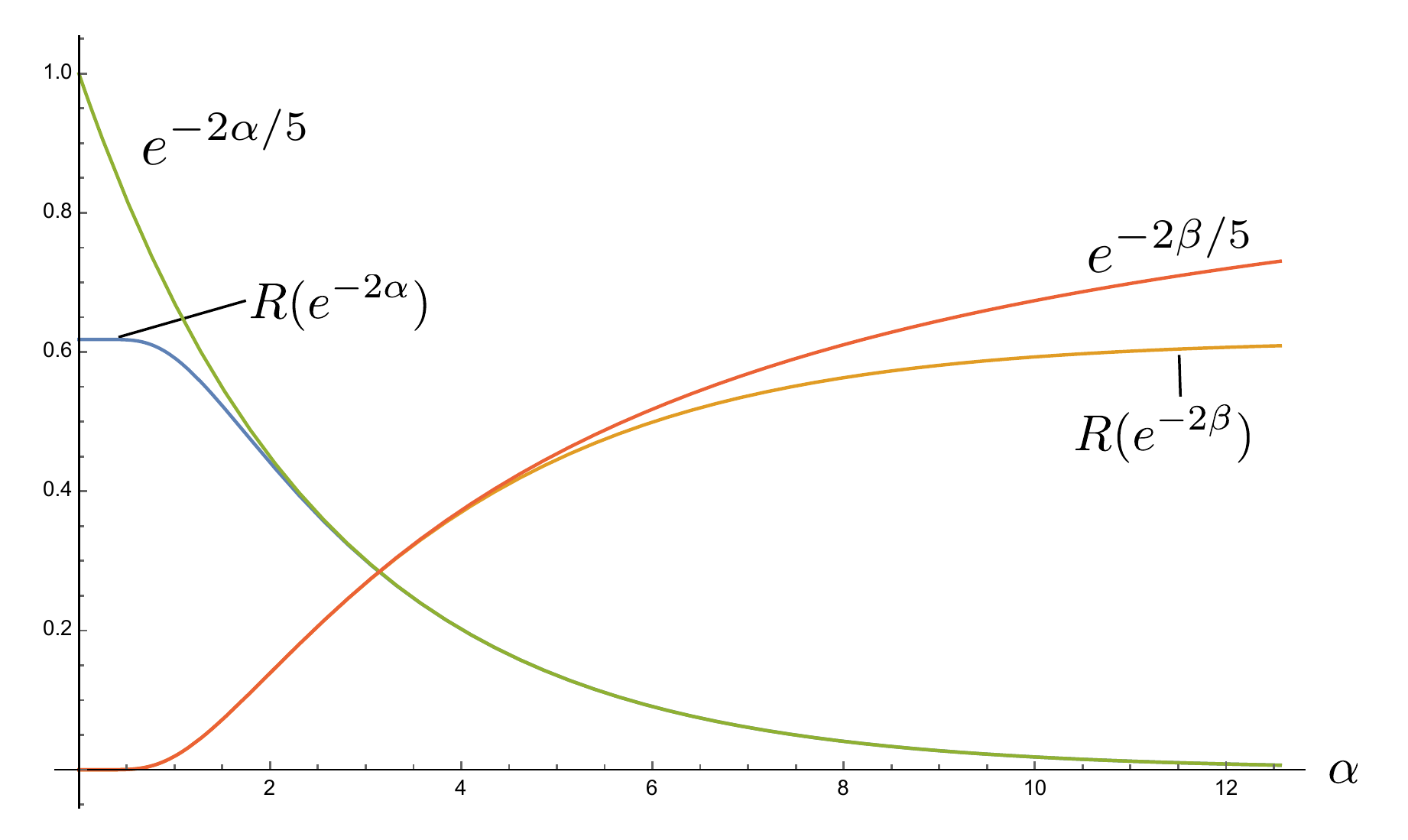}\caption{Comparisons of $e^{-2\alpha/5}$ and $R\left(e^{-2\alpha}\right)$,  and $e^{-2\beta/5}$ and $R(e^{-2\beta})$, as functions of $\alpha$, where $\alpha\beta=\pi^2$.}
\label{fig3}
\end{center}
\end{figure}
\begin{figure}[htbp]
\begin{center}
\includegraphics[width=0.75\textwidth]{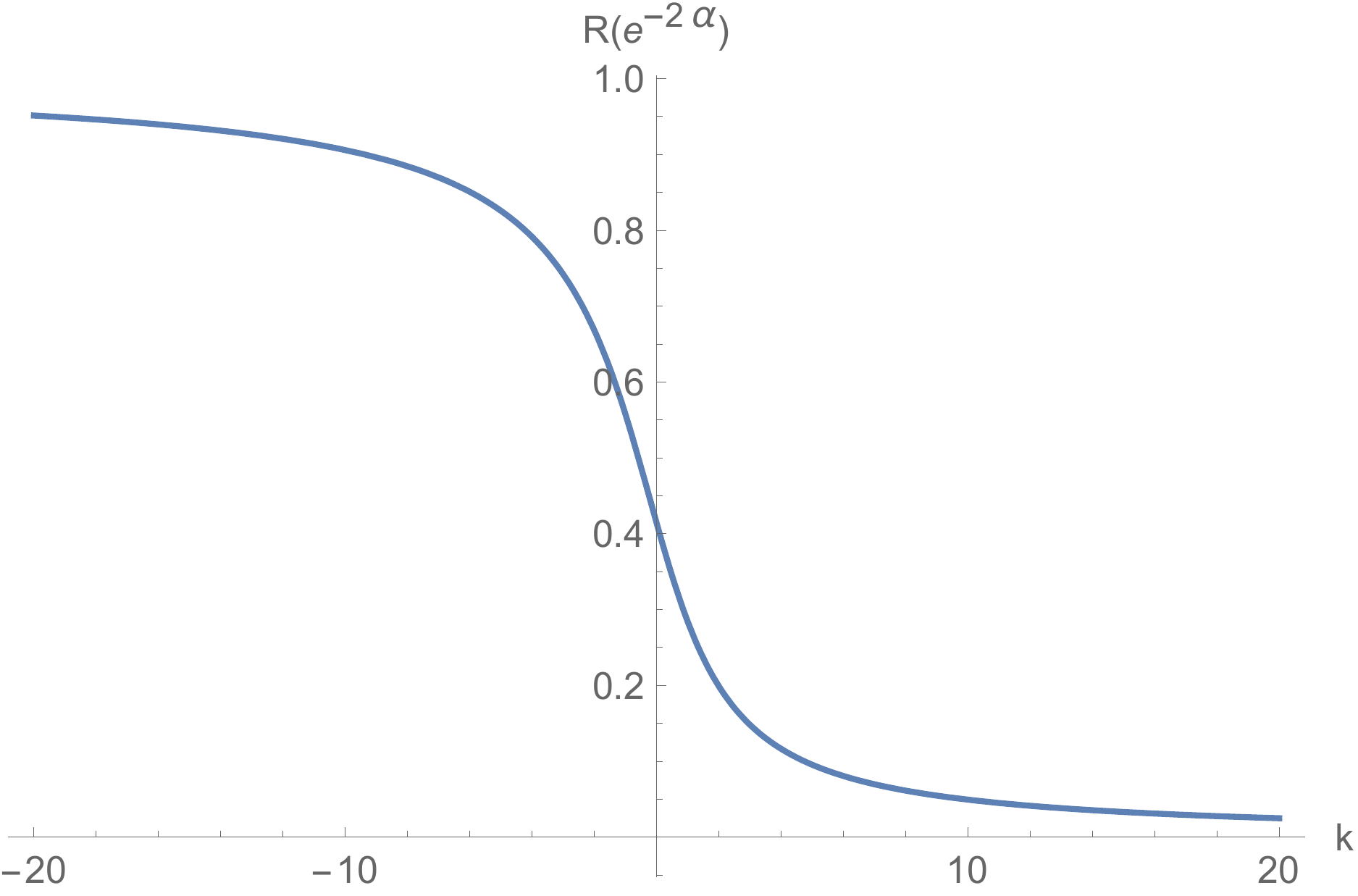}
\caption{$R\left(e^{-2\alpha}\right)$ as a function of $k$.}
\label{fig4}
\end{center}
\end{figure}

For example, consider $x_1=R\left(e^{-2\alpha}\right)$, $x_2=R(e^{-2\beta})$ and $2r=R\left(e^{-2\alpha}\right)+R(e^{-2\beta})$, where $\alpha\beta=\pi^2$. For each value of $r$, there are corresponding values of $\alpha$ and $\beta$ for which equation (28) has the roots $x_1=R\left(e^{-2\alpha}\right)$ and $x_2=R(e^{-2\beta})$. Figure 1 plots the values of $2r$ as a function of $\alpha$,  where $c=\phi$ and $\alpha\beta=\pi^2$. As the value of $\alpha$ increases from zero to $\pi$, the value of $2r$ decreases from $1/\phi$ to $2(\sqrt{\phi^2+1}-\phi)$. The minimum value of $2r$ is obtained when $\alpha=\beta=\pi$. 

Figure 2 shows how the values of $R\left(e^{-2\alpha}\right)$ and $R(e^{-2\beta})$ change with $\alpha$.  As $\alpha$ increases from zero to $\pi$, $R\left(e^{-2\alpha}\right)$ decreases from $1/\phi$ to $\sqrt{\phi^2+1}-\phi$ and $R(e^{-2\beta})$ increases from zero to $\sqrt{\phi^2+1}-\phi$.  The horizontal line in Figure 2 passes through $r=R(e^{-2\pi})$, which corresponds to $\alpha=\beta=\pi$. It is sufficient to know the value of $R\left(e^{-2\alpha}\right)$ over $0<\alpha\le \pi$; the corresponding value of $R(e^{-2\beta})$ can be obtained by the value of  the symmetric function around $r=R(e^{-2\pi})$, shown by the yellow line in Figure 2. Observe that $R\left(e^{-2\alpha}\right)$ and $R(e^{-2\beta})$ change very little when $0<\alpha<\pi/5$. Thereafter, the functions change almost linearly until $\alpha=\pi$. Figures \ref{fig3} compares the values of $R\left(e^{-2\alpha}\right)$ and  $e^{-2\alpha/5}$, and $R(e^{-2\beta})$ and  $e^{-2\beta/5}$, over the interval $0<\alpha\le 4\pi$.  For $\alpha\ge \pi$, $R\left(e^{-2\alpha}\right)$ is very close to $e^{-2\alpha/5}$ and becomes increasingly closer as $\alpha$ increases; and for $0<\alpha\le \pi$, $R(e^{-2\beta})$ is very close to $e^{-2\beta/5}$ and becomes increasingly closer as $\alpha$ decreases. Using the reciprocity formula in equation (6), we can approximate  $R\left(e^{-2\alpha}\right)$ by
$$R\left(e^{-2\alpha}\right)\approx 
\begin{cases}
{\left(1-\phi e^{-{2\pi^2/(5\alpha)}}\right)\Big/\left(\phi+e^{-{2\pi^2/(5\alpha)}}\right)}&{\rm if}\ 0<\alpha\le \pi,\\
e^{{-2\alpha/5}}&{\rm if}\ \alpha\ge \pi.\\
\end{cases}$$
The maximum difference between the values of $R\left(e^{-2\alpha}\right)$ and this approximate value is less than 0.000531 and occurs when $\alpha=\pi$. Figure \ref{fig4} shows the relation between $R\left(e^{-2\alpha}\right)$ and $k$.   Similar plots and results can be obtained about the terms in the other pairs of reciprocity formulas.

\subsection{Explicit expressions for $R(q)$}
Using equations (6) and (8), we can write the reciprocity formulas for $R\left(e^{-2\alpha}\right)$ and $R^5\left(e^{-2\alpha}\right)$ as
$$R\left(e^{-2\alpha}\right)=\left({\phi+{1\over \phi}\over \phi+R\left(e^{-2\beta}\right)}-1\right)\phi={1-\phi R\left(e^{-2\beta}\right)\over \phi+R\left(e^{-2\beta}\right)},\ \alpha\beta=\pi^2\eqno(29)\qquad\ \ $$
and
$$R^5\left(e^{-2\alpha}\right)=\left({\phi^5+{1\over \phi^5}\over \phi^5+R^5\left(e^{-2\beta}\right)}-1\right)\phi^5={1-\phi^5 R^5\left(e^{-2\beta}\right)\over \phi^5+R^5\left(e^{-2\beta}\right)},\  
\alpha\beta={\pi^2\over 5}.\eqno(30)$$
Starting with a pair of values for $R\left(e^{-2\alpha}\right)$ and $R(e^{-2\beta})$, where $\alpha\beta=\pi^2$,  we can iterate between equations (29) and (30) to obtain a sequence of values for $R(e^{-2\cdot 5^n \alpha})$, $R(e^{-2\cdot 5^n \beta})$, $R(e^{-2\pi^2/(5^n \alpha)})$
and $R(e^{-2\pi^2/(5^n \beta)})$, for all $n\ge 1$. These values have notable patterns when expressed in terms of $\phi$. Below, we obtain these expressions for 
$$R\left(e^{-2\pi(\sqrt{5})^n}\right)={1-\phi R\left(e^{-2\pi/(\sqrt{5})^n}\right)\over \phi+R\left(e^{-2\pi/(\sqrt{5})^n}\right)}$$
and
$$R^5\left(e^{-2\pi/(\sqrt{5})^{n+2}}\right)
={1-\phi^5 R^5\left(e^{-2\pi(\sqrt{5})^{n}}\right)\over \phi^5+R^5\left(e^{-2\pi(\sqrt{5})^{n}}\right)}.$$
Let
$$R\left(e^{-2\pi/(\sqrt{5})^n}\right)={A_n\over B_n}.$$
Let  $C_n=B_n-\phi A_n$ and $D_n=\phi + B_n$. Then
$$R\left(e^{-2\pi(\sqrt{5})^n}\right)={C_n\over D_n}$$
and
$$R^5\left(e^{-2\pi/(\sqrt{5})^{n+2}}\right)={D_n^5-\left(\phi C_n\right)^5\over \left(\phi D_n\right)^5+C_n^5}.$$

\bigskip\noindent
(1) Consider
$$R\left(e^{-2\pi}\right)={\sqrt{\phi^2+1}-\phi},$$
which is obtained by substituting $x_1=x_2=R(e^{-2\pi})$ in equation (6). Then
$$R^5\left(e^{-2\pi/5}\right)={{1-\phi^5 R^5\left(e^{-2\pi}\right)\over \phi^5+R^5\left(e^{-2\pi}\right)}}\qquad\qquad\qquad\qquad\quad$$
$$={{1-\phi^5 \left({\sqrt{\phi^2+1}-\phi}\right)^5\over \phi^5+\left({\sqrt{\phi^2+1}-\phi}\right)^5}}$$
and
$$R\left(e^{-2\cdot 5\pi}\right)={1-\phi R\left(e^{-2\pi/5}\right)\over \phi+R\left(e^{-2\pi/5}\right)}\qquad\qquad\qquad\qquad\qquad\qquad\qquad\qquad\qquad\qquad\quad$$
$$\ \ \ ={\sqrt[5]{{\phi^5+\left({\sqrt{\phi^2+1}-\phi}\right)^5}}-\phi \sqrt[5]{1-\phi^5 ({\sqrt{\phi^2+1}-\phi})^5}\over \phi \sqrt[5]{{\phi^5+\left({\sqrt{\phi^2+1}-\phi}\right)^5}}+\sqrt[5]{1-\phi^5 ({\sqrt{\phi^2+1}-\phi})^5}}\ .$$

\bigskip\noindent
(2) Consider
$$R^5\left(e^{-2\pi/\sqrt{5}}\right)={{\sqrt{\phi^{10}+1}-\phi^5}},$$
which is obtained by substituting $x_1=x_2=R^5(e^{-2\pi/\sqrt{5}})$ in equation (8). Then
$$R\left(e^{-2\pi\sqrt{5}}\right)
={1-\phi R\left(e^{-2\pi/\sqrt{5}}\right)\over \phi+R\left(e^{-2\pi/\sqrt{5}}\right)}\qquad\qquad\qquad\qquad\qquad\qquad\qquad\qquad\qquad$$
$$={1-\phi \sqrt[5]{{\sqrt{\phi^{10}+1}-\phi^5}}\over \phi +\sqrt[5]{{\sqrt{\phi^{10}+1}-\phi^5}}},\qquad\qquad\qquad\qquad\qquad$$
which Ramanujan expressed as follows in his second letter to Hardy (\cite{H}, p. 143):
$${1\over 1+}{e^{-2\pi\sqrt{5}}\over 1+}{e^{-4\pi\sqrt{5}}\over 1+ \dots}=\left[{\sqrt{5}\over 1+\sqrt[5]{5^{3/4}\left({\sqrt{5}-1\over 2}\right)^{5/2}-1}}-{\sqrt{5}+1\over 2}\right]e^{2\pi/\sqrt{5}}.$$
We use the expression for $R\left(e^{-2\pi\sqrt{5}}\right)$ to obtain
$$R^5\left(e^{-2\pi\over 5\sqrt{5}}\right)
={1-\phi^5 R^5\left(e^{-2\pi\sqrt{5}}\right)\over \phi^5+R^5\left(e^{-2\pi\sqrt{5}}\right)}\qquad\qquad\qquad\qquad\qquad\qquad\qquad\qquad\quad$$
$$\qquad \qquad\quad ={{\left(\phi + \sqrt[5]{{\sqrt{\phi^{10}+1}-\phi^5}}\right)^5-\phi^5 \left(1-\phi \sqrt[5]{{\sqrt{\phi^{10}+1}-\phi^5}}\right)^5\over \phi^5\left(\phi + \sqrt[5]{{\sqrt{\phi^{10}+1}-\phi^5}}\right)^5+\left(1-\phi \sqrt[5]{{\sqrt{\phi^{10}+1}-\phi^5}}\right)^5}}$$
and then use this expression to obtain
$$R\left(e^{-2\cdot 5\sqrt{5}\pi}\right)
={1-\phi R\left(e^{-2\pi/(5\sqrt{5})}\right)\over \phi+R\left(e^{-2\pi/(5\sqrt{5})}\right)}\qquad\qquad\qquad\qquad\qquad\qquad\qquad\qquad\qquad\qquad\qquad\qquad\qquad\qquad$$
$$\tiny ={\sqrt[5]{\phi^5 \left(\phi + \sqrt[5]{\sqrt{\phi^{10}+1}-\phi^5}\right)^5+
\left(1-\phi \sqrt[5]{\sqrt{\phi^{10}+1}-\phi^5}\right)^5}
-\phi \sqrt[5]{{\left(\phi + \sqrt[5]{\sqrt{\phi^{10}+1}-\phi^5}\right)^5-
\phi^5\left(1-\phi \sqrt[5]{\sqrt{\phi^{10}+1}-\phi^5}\right)^5}} \over 
\phi\sqrt[5]{\phi^5 \left(\phi + \sqrt[5]{\sqrt{\phi^{10}+1}-\phi^5}\right)^5+
\left(1-\phi \sqrt[5]{\sqrt{\phi^{10}+1}-\phi^5}\right)^5}
+ \sqrt[5]{{\left(\phi + \sqrt[5]{\sqrt{\phi^{10}+1}-\phi^5}\right)^5-
\phi^5\left(1-\phi \sqrt[5]{\sqrt{\phi^{10}+1}-\phi^5}\right)^5}} }\ .$$

\bigskip\noindent
(3) Consider
 $$R(e^{-\pi})=\frac{-\sqrt{\sqrt{5}} \left(\sqrt{5}+5\right)+\sqrt{2 \sqrt{5}+5}+\sqrt{5 \left(2 \sqrt{5}+5\right)}}{2 \left(\sqrt{5}-\sqrt{2 \sqrt{5}+5}\right)},$$
 which can be expressed as follows in terms of $\phi$:
$$R(e^{-\pi})={\phi+{1\over \phi}-\phi\sqrt{\phi}\over \sqrt{\phi}-{1\over \phi}\sqrt{\phi+{1\over \phi}}}.$$
Starting with this expression, we can iterate between equations (29) and (30) to obtain the values of  $R(e^{-4\cdot 5^{n-1}\pi})$,  $R(e^{-4\pi/5^n})$ and $R(e^{-\pi/5^n})$, for $n\ge 1$. The first such values are
$$R(e^{-4\pi})={1-\phi R(e^{-\pi})\over \phi +R(e^{-\pi})}\qquad\qquad\qquad\qquad\qquad\qquad\qquad\qquad\qquad$$
$$\ \ ={\left(\sqrt{\phi}-{1\over \phi}\sqrt{\phi + {1\over \phi}}\right) -\phi\left(\phi + {1\over \phi}-\phi\sqrt{\phi}\right)\over \phi+{1\over \phi}-\sqrt{\phi+{1\over \phi}}}\qquad\quad$$
$$R^5\left(e^{-4\pi/5}\right)={{1-\phi^5 R^5\left(e^{-\pi}\right)\over \phi^5+R^5\left(e^{-\pi}\right)}}\qquad\qquad\qquad\qquad\qquad\qquad\qquad\qquad\qquad$$
$$={\left(\sqrt{\phi}-{1\over \phi}\sqrt{\phi + {1\over \phi}}\right)^5 -\phi^5\left(\phi + {1\over \phi}-\phi\sqrt{\phi}\right)^5\over\phi^5\left(\sqrt{\phi}-{1\over \phi}\sqrt{\phi + {1\over \phi}}\right)^5+\left(\phi + {1\over \phi}-\phi\sqrt{\phi}\right)^5}\ \ $$
and
$$R^5\left(e^{-\pi/5}\right)={{1-\phi^5 R^5\left(e^{-4\pi}\right)\over \phi^5+R^5\left(e^{-4\pi}\right)}}\qquad\qquad\qquad\qquad\qquad\qquad\qquad\qquad\qquad\qquad$$
$$\qquad\qquad\qquad\qquad\  ={\left(\phi+{1\over \phi}-\sqrt{\phi+{1\over \phi}}\right)^5-\phi^5\left(\left(\sqrt{\phi}-{1\over \phi}\sqrt{\phi + {1\over \phi}}\right) -\phi\left(\phi + {1\over \phi}-\phi\sqrt{\phi}\right)\right)^5
\over 
\phi^5 \left(\phi+{1\over \phi}-\sqrt{\phi+{1\over \phi}}\right)^5+\left(\left(\sqrt{\phi}-{1\over \phi}\sqrt{\phi + {1\over \phi}}\right) -\phi\left(\phi + {1\over \phi}-\phi\sqrt{\phi}\right)\right)^5}\ .$$

\section{Ramanujan's cubic continued fraction and the Ramanujan-Selberg continued fraction}
The preceding analysis extends to the following continued fraction, which appears in entry 3.3.1 of Ramanujan's lost notebook (\cite{rln}, p. 94) and is now known as Ramanujan's cubic continued fraction:
$$G(q)={q^{1/3}\over 1+}\  {q+q^2\over 1+}\ {q^2+q^4\over 1+}\ {q^4+q^6\over 1+}\dots$$
$$\qquad\qquad\qquad\quad\ \  =q^{1/3}\prod_{j=1}^\infty {(1-q^{6j-5})(1-q^{6j-1})
\over (1-q^{6j-3})^2},\ |q|<1.\quad\ \ \ $$

Consider the equation
$${2^k-x_1x_2\over x_1+x_2}=1,\eqno(31)$$
which can be re-written as
$$(1+x_1)(1+x_2)=1+2^k.\eqno(32)$$
Setting $k=\pm 1$ and $k=\pm 3$ in equation (32) yields the following four reciprocity formulas for Ramanujan's cubic continued fraction. Table \ref{tab2} shows the terms associated with $x_1$ and $x_2$ in a reciprocity formula for each value of $k$.
\begin{table}[htp]
\caption{Values of $k$ associated with Ramanujan's cubic continued fraction}
\begin{center}
\begin{tabular}{rrrr}
\hline
$k$&$x_1$&$x_2$&$\alpha\beta$\\
\hline
$1$&$1/G\left(-e^{-\alpha}\right)$&$1/G(-e^{-\beta})$&$\pi^2$\\
$-1$&$G(e^{-\alpha})$&$G(e^{-\beta})$&$2\pi^2$\\
$3$&$1/G^3\left(-e^{-\alpha}\right)$&$1/G^3(-e^{-\beta})$&$\pi^2/3$\\
$-3$&$G^3(e^{-\alpha})$&$G^3(e^{-\beta})$&$2\pi^2/3$\\
\hline
\end{tabular}
\end{center}
\label{tab2}
\end{table}

\bigskip\noindent
1. Let $k=1$. Adiga, Kim, Naika and Madhusudhan \cite{aknm} obtained the following reciprocity formula in which $x_1=1/G(-e^{-\alpha})$ and  $x_2=1/G(-e^{-\beta})$:
$$\left(1+{1\over G(-e^{-\alpha})}\right)\left(1+{1\over G(-e^{-\beta})}\right)=1+2,\ \ \alpha\beta=\pi^2.\eqno(33)$$

\bigskip\noindent
2. Let $k=-1$. In his lost notebook, Ramanujan (\cite{rln}, Theorem 3.3.2) recorded the following reciprocity formula in which $x_1=G(e^{-\alpha})$ and  $x_2=G(e^{-\beta})$:
$$\Big(1+G(e^{-\alpha})\Big)\Big(1+G(e^{-\beta})\Big)=1+{1\over 2},\ \ \alpha\beta=2\pi^2.\eqno(34)$$

\bigskip\noindent
3. Let $k=3$. Adiga, Kim, Naika and Madhusudhan \cite{aknm} obtained the following reciprocity formula in which $x_1=1/G^3(-e^{-\alpha})$ and $x_2=1/G^3(-e^{-\beta})$:
$$\left(1+{1\over G^3(-e^{-\alpha})}\right)\left(1+{1\over G^3(-e^{-\beta})}\right)=1+2^3,\ \ \alpha\beta={\pi^2\over 3}.\eqno(35)$$

\bigskip\noindent
4. Let $k=-3$. Cooper(\cite{coop}, Theorem 6.12, p. 373)  obtained the following reciprocity formula in which $x_1=G^3(e^{-\alpha})$ and  $x_2=G^3(e^{-\beta})$:
$$\Big(1+G^3(e^{-\alpha})\Big)\Big(1+G^3(e^{-\beta})\Big)=1+{1\over 2^3},\ \ \alpha\beta={2\pi^2\over 3}.\eqno(36)$$

\bigskip
Let $x_1+x_2=2r$, where $x_1$ and $x_2$ are a pair of values in Table \ref{tab2}.  For each value of $k$ in the table, $x_1$ and $x_2$ are the roots $r\pm \sqrt{r^2+2r-2^k}$ of the quadratic equation  $x^2-2rx-(2r-2^k)=0$. Let $r\le 2^{k-1}$ so that both the roots are non-negative. If one root is zero, the other obtains it's largest possible value of $2^k$. The two roots are equal if $r^2+2r - 2^k=0$, in which case  $x_1=x_2=r=\sqrt{2^k+1}-1$.  If the two roots are different, the larger has a value between $2^k$ and $\sqrt{2^k+1}-1$, and the smaller a value between $\sqrt{2^k+1}-1$ and zero. Starting with a pair of values for $G\left(e^{-\alpha}\right)$ and $G(e^{-\beta})$, we can use equations (33) and (35), and equations  (34) and (36), to iteratively obtain sequences of values for $G(q)$. 

\bigskip
Finally, consider the following continued fraction, which Ramanujan recorded in his second notebook (\cite{Ram}, Chapter XIX, p. 225):
$$s(q)={q^{1/8}\over 1+}\  {q\over 1+q+}\ {q^2\over 1+q^2+}\ {q^3\over 1+q^3+}\dots, |q|<1.$$
 This continued fraction is now called the Ramanujan-Selberg continued fraction. Saikia (\cite{S}, Theorem 5.3, p. 52) obtained a reciprocity formula that can be expressed as
$$s^8(e^{-\alpha})+s^8(e^{-\beta})={1\over 16},\ \alpha\beta=\pi^2.$$
Thus, unlike the other preceding continued fractions, the Ramanujan-Selberg continued fraction is a solution to a linear equation, not a quadratic equation. Any positive number between zero and $1/16$ can be expressed as $x=s^8\left(e^{-\alpha}\right)$, and then ${1\over 16}-x$ can be expressed  as $s^8(e^{-\pi^2/\alpha})$.

\end{document}